\documentclass[10pt]{article}
\usepackage{a4}
\title{Unimodular lattices in dimensions 14 and 15 over the Eisenstein integers}
\author{Kanat Abdukhalikov\thanks{Supported by the Alexander von Humboldt Foundation}\\
{\normalsize Institute of Mathematics}\\
{\normalsize Pushkin Str 125}\\
{\normalsize Almaty 050010, Kazakhstan}\\
{\normalsize \texttt{abdukhalikov@math.kz}}
\and Rudolf Scharlau\\
{\normalsize Department of Mathematics}\\
{\normalsize University of Dortmund}\\
{\normalsize 44221 Dortmund, Germany}\\
{\normalsize \texttt{Rudolf.Scharlau@math.uni-dortmund.de}} }

\date{ }

\usepackage{amssymb}

\begin{document}

\newcommand{\F}{\mbox{$\mathbb{F}_4$}}          
\newcommand{\Z}{\mbox{$\mathbb{Z}$}}            
\newcommand{\Zo}{\mbox{$\mathbb{Z}[\omega ]$}}  
\newcommand{\Q}{\mbox{$\mathbb{Q}$}}            
\newcommand{\la}{\langle}
\newcommand{\ra}{\rangle}

\newtheorem{theorem}{Theorem}
\newtheorem{lemma}{Lemma}
\newtheorem{propo}{Proposition}

\maketitle

\begin{abstract}
All indecomposable unimodular hermitian lattices in dimensions 14
and 15 over the ring of integers in $\Q(\sqrt{-3})$ are determined.
Precisely one lattice in dimension 14 and two lattices in dimension
15 have minimal norm 3.
\end{abstract}

Keywords: Integral lattice, hermitian lattice, extremal lattice,
unimodular lattice, root system.

\bigskip\noindent
In 1978 W. Feit \cite{Feit} has classified the unimodular hermitian lattices
of dimensions up to 12 over the ring $\Zo$ of Eisenstein integers, where
$\omega$ is primitive third root of unity.
These lattices all have roots, that is, vectors of norm 2.
In dimension 13, for the first time a unimodular lattice without roots
appears \cite{Ab,Bach}. In \cite{Ab2} the unimodular lattices in dimension 13
are completely classified.  The root-free lattice turns out to be unique.  It
has minimal norm 3, and its automorphism group is isomorphic to the group
${\Z}_6\times {\rm PSp}_6(3)$ of order
$2^{10}\cdot 3^{10}\cdot 5\cdot 7\cdot 13$.
The remaining lattices all have roots; the rank of the root system is 12
in all cases.

In this paper, we classify the unimodular lattices in dimensions 14 and
15.   There are exactly 58, respectively 259 classes of indecomposable lattices
in these dimensions.
Below, we list their root systems and the orders of their automorphism
groups. Gram matrices for all lattices are available electronically via \\
\texttt{www.mathematik.uni-dortmund.de/$\sim$scharlau} .
There is only one root-free unimodular lattice of rank 14,
and there are two root-free unimodular lattices of rank 15.

The  lattices without roots have minimal norm 3; they are {\em extremal} as
introduced for unimodular Eisenstein lattices in \cite{Con}, Chapter 10.7.
They give rise to 3-modular extremal $\Z$-lattices in twice the dimension, as
defined by Quebbemann in \cite{Que}. See \cite{Con,Sch,Sch2} for more
information on extremal and modular lattices and their relation to modular
forms. In this context, the lattices classified in this paper can be
considered as complex structures on (extremal) 3-modular lattices. The
question for existence, uniqueness, and possibly a full classification of
extremal modular lattices has
been an ongoing challenge, both computationally and theoretically, after the
appearance of the influential paper \cite{Que}.

\medskip
Let $V$ be a vector space over $\Q(\sqrt{-3})$ with a positive definite
hermitian product $(,)$.  A {\em lattice} $L$ in $V$ is a finitely generated
$\Zo$-module contained in $V$ such that $L$ contains a basis of $V$ and
$(x,y)\in\Zo$ for all $x$, $y\in L$.  More precisely, one calls this an {\em
integral lattice}.  The ring $\Zo$ is a principal ideal domain. Thus every
finitely generated torsion-free $\Zo$-module is free.  The {\em discriminant}
$d(L)$ of $L$ is the determinant of the Gram matrix $(a_i,a_j)$ with respect
to some basis $a_1, \dots,a_n$ of $L$.
The lattice is {\em unimodular} if $d(L)=1$.  The norm of a vector $x\in L$ is
$N(x)=(x,x)$. The minimal norm of a lattice $L$ is $\min \{N(x,x) \mid x\in L,
\ x\neq 0 \}$.  Let $L^{(2)}$ denote the sublattice of $L$ generated by all
vectors of norm 2 in $L$, the so-called {\em roots} of $L$, and
let $\mu_2(L)$ denote the number of roots in $L$.  Finally $G(L)$ denotes the
group of all automorphisms of $L$ which preserve the form. The group $G(L)$ is
finite.  Every lattice can be uniquely decomposed into an orthogonal sum of
orthogonally indecomposable lattices. In view of the known classification in
smaller dimensions it is therefore sufficient to list the indecomposable
lattices.

\medskip We now come to the actual construction of unimodular
lattices. We use various methods, such as representations of finite
groups, lifting of linear codes, neighbour steps at suitable primes,
and computer-assisted as well as hand computations of various kinds.
We use extensively Schiemann's computer program \cite{Schi} for
calculation of invariants of hermitian lattices and their neighbours,
and the Magma Computational Algebra System \cite{Mag}. The crucial
point is the computation of automorphism groups, which (in both systems)
is based on the work \cite{PlSou}.  After we have constructed (and
distinguished) sufficiently many lattices, we use the mass formula as
worked out in \cite{Feit} to check the completeness of our list.  An
easy to calculate but important invariant of our lattices are their
root systems. Due to the well known classification of complex
reflection groups \cite{Shep,Cohen}, there exist a few additional root
systems over $\Zo$ which are not defined over $\Z$.  Of course, they
are familiar from the work on lattices of smaller dimension cited
above.  We use the following standard notation for root lattices.
$I_n$ denotes a lattice of rank $n$ with an orthonormal basis.
Equivalently
      $$I_n = \la (a_1, \dots, a_n) \ | \ a_i\in \Zo \ra$$
and if $x=(a_1, \dots, a_n)$, $y=(b_1, \dots, b_n)$ then
$(x,y)=\sum_{i=1}^n a_i {\bar b_i}$. The lattices
 $A_n \subseteq I_{n+1}$,  $E_8$, $E_7$ and $E_6$ are as usual.
For $\alpha \in \Zo$ define $D_n(\alpha) \subseteq I_n$ by
    $$D_n(\alpha)=\la (a_1, \dots, a_n) \ | \
        a_i\in \Zo , \ \sum_{i=1}^n a_i \equiv 0 \pmod{\alpha} \ra.$$
Finally,
  $$U_5 = \la A_5, \ (1/\sqrt{-3})(1,\omega,\omega^2,1,\omega,\omega^2)\ra ,$$
  $$U_6 = \la D_6(\sqrt{-3}), \ (1/\sqrt{-3})(1,\dots,1)\ra .$$
  The lattice $U_6$ is unimodular, and $U_5$ has discriminant 2.

The root systems of the lattices classified in this paper are given in
Tables 1 - 3 below. For the most of lattices in dimension 14, we give ``glue
vectors'' (additional generators) \cite{Con} together with the root systems.

In Table~1 we have collected separately all lattices
which are constructed from self-dual codes of length 14 over $\F$.
These codes have been classified in \cite{Con2,Mac}.
Let $\varphi$ denote the composite mapping $\Zo \rightarrow \Zo/2\Zo \cong \F$.
Then for any self-dual code $C$ of length 14 over $\F$ the lattice
$$
L=\frac{1}{\sqrt{2}}\la (a_1,\dots,a_{14})\in I_{14} \ |
                  \ (\varphi(a_1),\dots,\varphi(a_{14}))\in C \ra
$$ is a unimodular lattice of rank 14.  In particular, lattice No. 10 comes
from a quadratic residue code \cite{Lint,Ward} of length 14.

\medskip

Now we come to a particular construction of a unimodular lattice of rank 15.
\begin{lemma}
\label{lem}
The exterior square $U_6\wedge U_6$ is a unimodular lattice of rank
$15$ and minimal norm $3$.
\end{lemma}
\textsc{Proof:} Our proof is similar to the proof of Theorem 2.1 from
\cite{Cou}. Any element $a\in U_6\wedge U_6$ can be written in the
form $a=\sum^r_{i=1}x_i\wedge y_i$, where $r\leq 3$ and
$x_1,\dots,x_r$, $y_1,\dots,y_r$ are $2r$ linearly independent
vectors in $U_6$. The lattice $L_{2r}$ generated by
$x_1,\dots,x_r$, $y_1,\dots,y_r$ depends only on $a$. We have
\begin{eqnarray*}
(a,a)
& = & \left( \sum^r_{i=1}x_i\wedge y_i, \sum^r_{j=1}x_j\wedge y_j\right) \\
& = & \sum_{i,j}[(x_i,x_j)(y_i,y_j) - (x_i,y_j)(y_i,x_j)] \\
& = & {\rm Tr}(A\overline{B} - C\overline{C}),
\end{eqnarray*}
where $A=((x_i,x_j))_{1\leq i,j\leq r}$ ,
$B=((y_i,y_j))_{1\leq i,j\leq r}$ and $C=((x_i,y_j))_{1\leq i,j\leq r}$.
Define
$$H = \left(\begin{array}{cc}
                        A & C \\
          ^t{\overline C} & B
            \end{array}
\right), \ J = \left(\begin{array}{cc}
                           0 & I_r \\
                        -I_r & 0
                   \end{array}
\right).
$$
Then $H$ is the hermitian matrix of the lattice $L_{2r}$ and
$$
(a,a) = \frac{1}{2}{\rm Tr}(H {^t\overline{J}} \overline{H}J)
        \geq \frac{2r}{2}|H|^{1/2r} \cdot |{^t\overline{J}}\overline{H}J|^{1/2r}
          = r|H|^{1/r},
$$
where we used the inequality ${\rm Tr}(MN) \geq m|M|^{1/m}|N|^{1/m}$ valid
for any hermitian positive definite matrices $M$, $N$ of degree $m$ (it can be
proved similarly to the proof of Lemma 7.1.3 in \cite{Kit}).  Now we show that
$(a,a)\geq 3$.  Indeed, for $r=3$ we have $(a,a)\geq 3\cdot |H|^{1/3}\geq 3$;
for $r=2$ we have $(a,a)\geq 2\cdot |H|^{1/2} > 2$; for $r=1$ we have
$(a,a)\geq 1\cdot |H|=d(L_{2r})\geq 3$.
Finally, it is easy to construct an element $a$ of norm 3:
take $a=x \wedge y$ with $(x,x)=(y,y)=2$ and $(x,y)=1$.   \hfill{$\Box$}

\begin{theorem}
There are $58$ hermitian indecomposable unimodular lattices over $\Zo$ of
dimension $14$. They are listed in Tables $1$ - $2$.
One of them has no roots; the remaining lattices have root systems of ranks
$6$, $8$, $10$, $11$, $12$, $13$ and $14$.
The root-free lattice $L_0$ has minimal norm $3$ and automorphism group
${\Z}_6\times {\rm G}_2(3).2$ of order
$2^8\cdot 3^7\cdot 7\cdot 13$.
\end{theorem}

\begin{theorem}
There are $259$ hermitian indecomposable unimodular lattices
over $\Zo$ of dimension $15$. They are listed in
Table~$3$. For any integer $r$
from $0$ to $15$ there is an indecomposable unimodular lattice of
dimension $15$ with root system of rank $r$. There are only two
root-free lattices. They have minimal norm $3$. One of the root-free
lattices is isometric to the exterior square $U_6 \wedge U_6$ and
has automorphism group $\Z_2 \times 3.{\rm U}_4(3).2$ of order
$2^9\cdot3^7\cdot 5\cdot 7$. The second root-free lattice has
automorphism group $\Z_2 \times (3_+^{1+2}\times 3_+^{1+2}).{\rm SL}_2(3).2$
of order $2^5\cdot 3^7$.
\end{theorem}
\textsc{Proof:}
Let $\mathfrak{G}_n$ be a complete set of representatives of the isomorphisms
classes of lattices of rank $n$ and discriminant 1.
Let $\mathfrak{G}_n^0$ be the set of all those lattices in $\mathfrak{G}_n$
which contain no element $x$ with $(x,x)=1$. Define
$$
M_n = \sum_{L\in {\mathfrak{G}_n}} \frac{1}{|G(L)|}, \quad
        M_n(2) = \sum_{L\in {\mathfrak{G}_n}} \frac{\mu_2(L)}{|G(L)|},
$$
$$Y_n = \sum_{L\in \mathfrak{G}_n^0} \frac{1}{|G(L)|}, \quad
        Y_n(2) = \sum_{L\in \mathfrak{G}_n^0} \frac{\mu_2(L)}{|G(L)|}.
$$

These quantities can be calculated by in principle well known formulas: the
mass formula for positive definite hermitian lattices, and a similar formula
for average representation numbers like $\mu_2$. In our situation these
formulas have been made explicit in \cite{Feit}, and one can calculate

$$
M_{14}=\frac{689532652191539}{2^{25}\cdot 3^{19}\cdot 5^3\cdot 11}, \quad
    M_{14}(2)=\frac{1722885336811913}{2^{23}\cdot 3^{17}\cdot 5^3\cdot 11},
$$
$$
Y_{14}=\frac{902121810728981}
               {2^{24}\cdot 3^{14}\cdot 5^3\cdot 7^2\cdot 11\cdot 13}, \quad
    Y_{14}(2)=\frac{321547203435163}
               {2^{22}\cdot 3^{12}\cdot 5^3\cdot 11\cdot 13}.
$$
$$
M_{15}=\frac{4366489808207046403}{2^{26}\cdot 3^{21}\cdot 5^3\cdot 11},\quad
    M_{15}(2)=\frac{1884491476714586441}{2^{24}\cdot 3^{18}\cdot 5^3\cdot 11},
$$
$$
Y_{15}=\frac{3619970721202760389}
                {2^{18}\cdot 3^{20}\cdot 5^3\cdot 7^2\cdot 11\cdot 13},\quad
    Y_{15}(2)=\frac{312324214206248801}
               {2^{16}\cdot 3^{17}\cdot 5^2\cdot 7^2\cdot 11\cdot 13}.
$$

Now we use an explicit list of Gram matrices of the
constructed lattices.
We use Schiemann's computer program and the Magma Computational
Algebra System to calculate $\mu_2(L)$ and $|G(L)|$ in each case.
Then a  straightforward check of the mass formulas for $Y_{14}$
and $Y_{15}$ confirms the completeness of our list. The
formulas for $Y_{14}(2)$ and $Y_{15}(2)$ give an additional check (ignoring the
lattices without roots).

We now determine the structure of the automorphism groups of the three
root-free lattices.  The group ${\rm G}_2(3).2$ has a representation of
dimension 14 over $\Q(\sqrt{-3})$ (see \cite{Atl}) and it acts on a lattice
$L$ of rank 14 over $\Zo$. Now results of \cite{Nebe} on rational maximal
finite matrix groups in dimension 28 imply that this lattice is unimodular,
and therefore is isometric to our root-free lattice $L_0$.
The factor $\Z_6$ in $G(L_0)$ comes from roots of unity in $\Zo$.

Now we come to the automorphism group of the root-free lattice $U_6 \wedge
U_6$ introduced in Lemma 1.
Recall that ${\rm Aut}(U_6) \cong 6.{\rm U}_4(3).2$.
There is a natural homomorphism
$\varphi : {\rm Aut}(U_6) \mapsto {\rm Aut}(U_6 \wedge U_6)$, given by
$\varphi(\sigma)(x\wedge y) = \sigma(x)\wedge \sigma(y)$ with
${\rm  Ker}\ \varphi=\{\pm {\rm id}\}$.
Therefore, ${\rm Aut}(U_6 \wedge U_6) \supseteq 3.{\rm U}_4(3).2$.
But $\{\pm {\rm id}\}\cdot \varphi({\rm Aut}(U_6))$ is a group of
order $2^9\cdot3^7\cdot 5\cdot 7$. The known orders of the two automorphism
groups of root-free lattices in dimension 15 show that this group must be the
full group.
Therefore, ${\rm Aut}(U_6 \wedge U_6) \cong \Z_2 \times 3.{\rm U}_4(3).2$.

The structure of the automorphism group of the second root-free lattice, of
order $2^5 \cdot 3^7$ was
determined with the help of Magma.  \hfill{$\Box$}

The extremal 3-modular $\Z$-lattice in dimension 28 coming from our root-free
lattice in dimension 14 appeared for the first time in \cite{Nebe} since its
automorphism group is a rational irreducible maximal finite subgroup of ${\rm
GL}_{28}(\Q )$. The two root-free lattices of rank 15 have been found  by
A. Schiemann (\cite{Schi}, also mentioned in \cite{Sch}), using a
sophisticated computer search for lattices with large minimum.

\medskip
We want to add some comments on algebraic techniques for
lattices which are related in one or the other way to the automorphism
group. One of the most important algebraic techniques for lattices is
``gluing'', which in the most general case means: (re)construction a lattice
$L$ from a known, ``simpler'' sublattice $M$. Then $L$ is generated by $M$ and
certain additional vectors, called glue vectors.  If $M$ has the same rank as
$L$, that is, generates the vector space $V$, the glue vectors must come from
the dual lattice $M^\#$, which can be analyzed once for all. Typically, $M$ is
decomposable while $L$ is not, so the additional vectors ``glue the components
of $M$ together'', this is where the terminology comes from. An important
choice for $M$ is the lattice $L^{(2)}$ generated by the roots. Then $M$ is
characteristic, that is $\mathrm{Aut}(L)$ acts on $M$, and is a subgroup of
$\mathrm{Aut}(M)$ if $M$ is of full rank. If $M$ is not of full rank, it must
be glued with a lattice on its orthogonal subspace. Here are two cases where
the rank of the root system is very small.

Let us consider the unimodular 15-dimensional lattice $\Lambda$ (No. 257) with
root system $A_2$. This lattice is obtained by gluing of a 2-dimensional
lattice with root system $A_2$ and a 13-dimensional lattice $\Lambda'$ with
discriminant 3. Computations on Magma and \cite{Atl} show that ${\rm
Aut}(\Lambda) \cong 6.({\rm S}_3\times {\rm PSL}_2(27).3)$. Therefore ${\rm
Aut}(\Lambda') \cong 6.({\rm PSL}_2(27).3)$ and $\Lambda'$ is associated with
the irreducible complex character of the group ${\rm PSL}_2(27)$ of degree
13. It is the Weil character and it can indeed be realized over $\Zo$ (see
\cite{Rie}).

Similar reasoning shows that the 15-dimensional unimodular lattice
No. 259 with root system $A_1$ has automorphism group $6.({\rm S}_2\times
{\rm PSL}_2(13))$ and it is obtained from a 14-dimensional lattice
of discriminant 2 with automorphism group $6.{\rm PSL}_2(13)$. This
lattice generates an irreducible complex character of the group
${\rm PSL}_2(13)$ of degree 14 (the case of complex character of degree
13 is excluded since the description \cite{Ab} of hermitian
lattices of dimension 13 associated with the irreducible complex
character of the group ${\rm PSL}_2(13)$ of degree 13 shows that
this case is impossible).

Of course, root sublattices and their complements are not the only choice for
characteristic sublattices. If the automorphism group is already known, one
could proceed as in the following example.
The complex representation associated with the action of the
automorphism group $\Z_2 \times (3_+^{1+2}\times 3_+^{1+2}).{\rm SL}_2(3).2$
on the 15-dimensional root-free lattice $L_2$ is reducible, and is a sum of
two irreducible representations of degrees 6 and 9 defined over
$\Q (\sqrt{-3})$.
And the complex
representation, associated with the action of the subgroup $\Z_2 \times
(3_+^{1+2}\times 3_+^{1+2}).{\rm SL}_2(3)$, is a sum of three irreducible
representations of degrees 3, 3 and 9.
This decomposition gives rise to invariant sublattices of $L_2$, and using the
general method of gluing it should be possible to obtain a computer-free
description of $L_2$.

Obviously, root lattices cannot be used directly to obtain lattices with
minimal norm 3 or larger. But one can use a scaled copy of a root lattice as
in the following example. We start with a root system over $\Z$, and in
addition to scaling we also use the larger ground ring $\Zo$.
Let $n$ be of the shape $n=2r\bar{r}$ for some  $r \in \Zo$.
Let $C$ be a self-dual code over $\F$ of length $n$  and
$$
\Lambda=\la (a_1,\dots,a_{n})\in A_{n-1} \ |
                  \ (\varphi(a_1),\dots,\varphi(a_{n}))\in C \ra
     + \frac{1}{r}\la (1-n,1,\dots,1)\ra.
$$
Then $\frac{1}{\sqrt{2}}\Lambda$ is a unimodular hermitian lattice of rank
$n-1$.
This construction had been used in \cite{Ab2} to obtain the root-free lattice
in dimension 13. It can also be applied in e.g. dimension 17.

\medskip
We close with some remarks on the theta series of our lattices.
There is one interesting observation in dimension 14 which
might deserve further attention.

\begin{propo}
Any  indecomposable unimodular lattice of rank $14$, more generally any
$3$-modular $\Z$-lattice of rank $28$ with minimum $\ge 2$ has
exactly $17472=2^6\cdot 3\cdot 7 \cdot 13$ vectors of norm $3$.
\end{propo}
\textsc{Proof:} (See \cite{Sch} for the method): The relevant algebra
$\mathcal M^*(3,\chi)$ of modular forms (where $\chi$ is an appropriate character) is
generated by the cusp form $\Delta_3(q)=\eta(q)^6 \, \eta(3q)^6$ of weight $6$ and
the theta-series $\theta_3:= \theta_{A_2}$ of the $A_2$-root lattice, of weight
$1$. A linear basis for the space $\mathcal M_{14}(3,\chi)$ of forms of weight
14, containing the
theta-series of the lattices in question is the following:
$$
\renewcommand{\arraystretch}{1.3}
\begin{array}{ll rrrrrrrr}
\theta_3^{14} &= &1 & + &84q &+ &3276q^2 &+ &78708q^3 &+ O(q^4)\\
{\theta_3}^8\, \Delta_3 &= &&& q &+ &42q^2 &+ &729q^3 &+ O(q^4)\\
{\theta_3}^2\, {\Delta_3} ^2 &= & & && & q^2 & & &+ O(q^4)
\end{array}
$$
The theta series of a lattice without vectors of norm 1 is of the form
$$ \theta_3^{14} - 84 {\theta_3}^8\, \Delta_3 + c {\theta_3}^2\, {\Delta_3} ^2
$$
with $c$ given by the number of vectors of norm 2. The absence of $q^3$ in the
third basis vector proves the claim. \hfill{$\Box$}

\medskip
Note that $M_{13}\approx  0.00000014$ (there are 14 lattices of rank 13),
$M_{14}\approx 0.000012$ (58 lattices), $M_{15}\approx 0.0045$ (259 lattices),
$M_{16}\approx 6.57$, $M_{17}\approx 42188.20$.
So beyond dimension 15, a full classification of unimodular Eisenstein
lattices might not be possible any more.

\renewcommand{\arraystretch}{1.2}
\renewcommand{\cdot}{.}

\begin{table}[t]{Table 1. Unimodular lattices of rank 14 (obtained from self-dual codes over $\F$)}
\label{table}
\begin{center}

\begin{tabular}{|l|l|l|l|l|}  \hline
No
& $L$
& $L^{(2)}$       &  $|G(L)|$  &  $\mu_2(L)$     \\ \hline
1
&$\la L^{(2)}, \frac{1}{2}(\sqrt{-3},1,\dots,1) \ra$
& $D_{14}(2)$
& $2^{24}. 3^6. 5^2. 7^2. 11. 13$ & 1092  \\ \hline
2
& $\la L^{(2)}, y_1 + y_2 \ra$,
& ${E_7}^2 $
& $2^{21}\cdot 3^9\cdot 5^2\cdot 7^2$ & 756   \\
& $y_i \in E_7^* - E_7$
&
&
&                                         \\ \hline
3
& $\la L^{(2)}, \frac{1}{2}(\sqrt{-3},1,\dots,1)+y_3$,
& $D_8(2)  U_5  A_1$
& $2^{22}\cdot 3^7\cdot 5^2\cdot 7$      & 612   \\
& $\frac{1}{2}(-\sqrt{-3},1,\dots,1)+y_2 \ra$,
&
&
&                                  \\
& $y_2=\frac{1}{2\sqrt{-3}}(-5,1,1,1,1,1)$, $y_3=\frac{1}{2}(-1,1)$
&
&
&                                         \\ \hline
4
& $\la L^{(2)},y_1+z,y_2+z \ra$
& ${U_5}^2  D_4(2)$
& $2^{21}\cdot 3^{11}\cdot 5^2$      & 612   \\
& $z=\frac{1}{2}(\sqrt{-3},1,1,1)$,
&
&
&                                  \\
& $y_i=\frac{1}{2\sqrt{-3}}(-5,1,1,1,1,1)$
&
&
&                                         \\ \hline
5
& $\la L^{(2)},$
& $D_8(2)  D_6(2)$
& $2^{23}\cdot 3^5\cdot 5^2\cdot 7$      & 516   \\
& $\frac{1}{2}(1,\dots,1) + (1,0,\dots,0)$,
&
&
&                                  \\
& $(1,0,\dots,0) + \frac{1}{2}(\sqrt{-3},1,\dots,1) \ra$
&
&
&                                         \\ \hline
6
& $\la L^{(2)},$
& ${D_6(2)}^2  {A_1}^2$
& $2^{21}\cdot 3^5\cdot 5^2$      & 372   \\
& $(1,0,\dots,0)+(1,0,\dots,0)+\omega y_1+\omega^2 y_2$,
&
&
&                                  \\
& $\frac{1}{2}(\sqrt{-3},1,\dots,1)_1+y_1+y_2$,
&
&
&                                         \\
& $\frac{1}{2}(\sqrt{-3},1,\dots,1)_1-\frac{1}{2}(\sqrt{-3},1,\dots,1)_2\ra,$
&
&
&                                         \\
& $y_i=\frac{1}{2}(-1,1)$
&
&
&                                         \\ \hline
7
& $\la L^{(2)},$
& $D_6(2)  {D_4(2)}^2 $
& $2^{22}\cdot 3^6\cdot 5$      & 324   \\
& $(1,0,\dots,0) + (1,0,0,0) + (1,0,0,0)$,
&
&
&                                  \\
& $(0,\dots,0,\omega) + \frac{1}{2}(1,1,1,1)+\frac{1}{2}(1,1,1,1)$,
&
&
&                                  \\
& $\frac{1}{2}(1,\dots,1) + \frac{1}{2}(\sqrt{-3},1,1,1) \ra$
&
&
&                                         \\ \hline
8
& $\la L^{(2)},$
& $D_4(2)^3  {A_1}^2 $
& $2^{21}\cdot 3^5$      & 228   \\
& $(1,0,0,0)_1 + (1,0,0,0)_2 + (1,0,0,0)_3$,
&
&
&                                  \\
& $y_1 + y_2 + y_4 + y_5$,
&
&
&                                  \\
& $y_2 + y_3 + \omega^2 y_4+\omega^2 y_5$,
&
&
&                                         \\
& $(1,0,0,0)_2 + (\omega^2,0,0,0)_3 + \omega y_4 + \omega^2 y_5 \ra$,
&
&
&                                  \\
& $y_{1,2,3}=\frac{1}{2}(1,1,1,1)$, $y_{4,5}=\frac{1}{2}(-1,1)$
&
&
&                                  \\       \hline
9
& $\la L^{(2)},$
& $D_4(2)^2 {A_1}^6$
& $2^{20}\cdot 3^5$      & 180   \\
& $(1,0,0,0)_1+y_3+y_4+\omega y_5+\omega y_6$,
&
&
&                                  \\
& $\frac{1}{2}(1,1,1,1)_1+\omega^2 y_2+\omega^2 y_3+y_4+y_5$,
&
&
&                                  \\
& $(1,0,0,0)_2+y_2+\omega y_3+y_5+\omega y_6$,
&
&
&                                         \\
& $\frac{1}{2}(1,1,1,1)_2+y_2+y_3+\omega^2 y_4+\omega^2 y_5$,
&
&
&                                  \\
& $y_1+y_2+y_3+y_4+y_5+y_6\ra$,
&
&
&                                  \\
& $y_i=\frac{1}{2}(-1,1)$
&
&
&                                  \\       \hline
10
& $\la L^{(2)}, \alpha_1y_1+\cdots +\alpha_{14}y_{14} \ra$
& ${A_1}^{14}$
& $2^{16}\cdot 3^2\cdot 7 \cdot 13$      & 84   \\
& $y_i=\frac{1}{2}(-1,1)$,
&
&
&                                  \\
& $(\alpha_1,\dots,\alpha_{14})_{\bmod{2}} \in QR(14)$
&
&
&                                         \\ \hline
\end{tabular}
\end{center}
\end{table}

\begin{table}[t]{Table 2. Unimodular lattices of rank 14}
\label{table2}
\begin{center}
\begin{tabular}{|l|l|l|l|l|}  \hline
No
& $L$
& $L^{(2)}$       &  $|G(L)|$  &  $\mu_2(L)$     \\ \hline
11
& $L_0$
& $\varnothing$             & $2^8\cdot 3^7\cdot 7\cdot 13$ & 0  \\ \hline
12
& $\la L^{(2)}, \frac{1}{\sqrt{-3}}(1,\dots,1)+y_2\ra ,$
& $D_8(\sqrt{-3})  E_6$
& $2^{15}\cdot 3^{14}\cdot 5^2\cdot 7$      & 720   \\
& $y_2 \in E_6^*-E_6$
&
&
&                                         \\ \hline
13
& $\la L^{(2)}, (2+\sqrt{-3})y_1 + y_2 \ra$,
& $A_{13}  A_1$
& $2^{13}\cdot 3^6\cdot 5^2\cdot 7^2\cdot 11\cdot 13$ & 552   \\
& $y_1=\frac{1}{14}(-13,1,\dots,1)$, $y_2=\frac{1}{2}(-1,1)$ & & &
\\ \hline
14
& Complex conjugate of 13
&     &                    &          \\ \hline
15
& $\la L^{(2)},(1,0,0,0,0,0)+y_2+y_3,$
& $D_6(\sqrt{-3})  E_6 A_2$
& $2^{13}\cdot 3^{13}\cdot 5^2$      & 504   \\
& $\frac{1}{\sqrt{-3}}(1,1,1,1,1,1)-\sqrt{-3}y_3 \ra$,
&
&
&                                  \\
& $y_2\in E_6^*-E_6$, $y_3=\frac{1}{3}(-2,1,1)$
&
&
&                                         \\ \hline
16
& $\la L^{(2)},$
& $A_{11}  A_3 $
& $2^{14}\cdot 3^7\cdot 5^2\cdot 7\cdot 11$      & 432   \\
& $\frac{\sqrt{-3}}{12}(-11,1,\dots,1) + \frac{\sqrt{-3}}{4}(-3,1,1,1) \ra$,
&
&
&                           \\ \hline
17
& $\la L^{(2)},\sqrt{-3}y_1 + y_2 \ra$,
& $A_8 E_6$
& $2^{15}\cdot 3^9\cdot 5^2\cdot 7$ & 432   \\
& $y_1=\frac{1}{9}(-8,1,\dots,1)$, $y_2\in E_6^*-E_6$
&
&
&                                         \\ \hline
18
& $\la L^{(2)},y_1+y_2+y_3$, $y_1-y_2+y_4 \ra$
& ${D_5(\sqrt{-3})}^2  {A_2}^2$
& $2^{11}\cdot 3^{13}\cdot 5^2$      & 396   \\
&  $y_{1,2}=\frac{1}{\sqrt{-3}}(1,\dots,1)$, $y_{3,4}=\frac{1}{3}(-2,1,1)$
&
&
&                                         \\ \hline
19
& $\la L^{(2)} \oplus \la x\ra,$
& $A_8  D_5(\sqrt{-3})$
& $2^{11}\cdot 3^{10}\cdot 5^2\cdot 7$      & 396   \\
& $y_1+y_2+\frac{2}{9}x$, $(1,0,0,0,0)+\frac{\sqrt{-3}}{3}x \ra$,
&
&
&                                  \\
& $y_1=\frac{1}{9}(-8,1,\dots,1)$,
&
&
&                                         \\
& $y_2=\frac{1}{\sqrt{-3}}(1,1,1,1,1)$, $(x,x)=9$
&
&
&                                         \\ \hline
20
& $\la L^{(2)},(2+\sqrt{-3})y_1 + y_2 \ra$,
& ${A_7}^2$
& $2^{15}\cdot 3^5\cdot 5^2\cdot 7^2$ & 336   \\
& $y_i=\frac{1}{8}(-7,1,\dots,1)$
&
&
&                                         \\ \hline
21
& $\la L^{(2)}, y_1+2\sqrt{-3}y_2+y_3\ra$
& $A_9  A_4  A_1$
& $2^{13}\cdot 3^6\cdot 5^3\cdot 7$ & 336   \\
& $y_1=\frac{1}{10}(-9,1,\dots,1)$,
&
&
&                                  \\
& $y_2=\frac{1}{5}(-4,1,\dots,1)$, $y_3=\frac{1}{2}(-1,1)$
&
&
&                                         \\ \hline
22
& $\la L^{(2)} \oplus \la x\ra,$
& $A_8  D_5(2)$
& $2^{15}\cdot 3^6\cdot 5^2 \cdot 7$      & 336   \\
& $(2+\sqrt{-3})y_1+\sqrt{-3}y_2+\frac{1}{36}x \ra$,
&
&
&                                  \\
& $y_1=\frac{1}{9}(-8,1,\dots,1)$,
&
&
&                                  \\
& $y_2=\frac{1}{2}(1,\dots,1)$, $(x,x)=36$
&
&
&                                         \\ \hline
23
& Complex conjugate of 22
&     &                    &          \\ \hline
24
& $\la L^{(2)}\oplus \la x\ra$,
& $A_6  D_5(2)  A_2$
& $2^{13}\cdot 3^5\cdot 5^2 \cdot 7$      & 264   \\
& $\sqrt{-3}y_1+\sqrt{-3}y_2+y_3+\frac{1}{84}x\ra$,
&
&
&                                  \\
&  $y_1=\frac{1}{7}(-6,1,\dots,1)$, $y_2=\frac{1}{\sqrt{-3}}(1,\dots,1)$,
&
&
&                                         \\
&  $y_3=\frac{1}{3}(-2,1,1)$, $(x,x)=84$
&
&
&                                         \\ \hline
25
& $\la L^{(2)} \oplus \la x\ra$,
& $A_7  A_5  A_1$
& $2^{13}\cdot 3^5\cdot 5^2 \cdot 7$      & 264   \\
& $y_1+y_3+\frac{2+\sqrt{-3}}{8}x$, $y_2+y_3+\frac{2-\sqrt{-3}}{6}x\ra$,
&
&
&                                  \\
&  $y_1=\frac{1}{8}(-7,1,\dots,1)$, $y_2=\frac{1}{6}(-5,1,\dots,1)$,
&
&
&                                         \\
&  $y_3=\frac{1}{2}(-1,1)$, $(x,x)=24$
&
&
&                                         \\ \hline
26
& Complex conjugate of 25
&     &                    &          \\ \hline
\end{tabular}
\end{center}
\end{table}


\begin{table}[t]{Table 2. Unimodular lattices of rank 14 (continued)}
\label{table3}
\begin{center}

\begin{tabular}{|l|l|l|l|l|}  \hline
No
& $L$
& $L^{(2)}$       &  $|G(L)|$  &  $\mu_2(L)$     \\ \hline
27
& $\la L^{(2)} \oplus \la x\ra$,
& $A_7  D_4(\sqrt{-3})  {A_1}^2$
& $2^{13}\cdot 3^7\cdot 5 \cdot 7$      & 288   \\
& $y_1+y_3+\frac{2+\sqrt{-3}}{8}x, y_3+y_4+\frac{1}{2}x,
y_2+\frac{1}{3}x\ra$,
&
&
&                                  \\
&  $y_1=\frac{1}{8}(-7,1,\dots,1)$, $y_2=\frac{1}{\sqrt{-3}}(1,1,1,1)$,
&
&
&                                         \\
&  $y_{3,4}=\frac{1}{2}(-1,1)$, $(x,x)=24$
&
&
&                                         \\ \hline
28
& $\la L^{(2)},$
& ${A_5}^2   D_4(\sqrt{-3})$
& $2^{13}\cdot 3^9\cdot 5^2$      & 288   \\
& $y_1+y_2+y_3, 2\sqrt{-3}(y_1-y_2) \ra,$
&
&
&                                  \\
& $y_{1,2}=\frac{1}{6}(-5,1,\dots,1)$, $y_3=\frac{1}{\sqrt{-3}}(1,1,1,1)$
&
&
&                                         \\ \hline
29
& $\la L^{(2)} \oplus \la x\ra,$
& $D_5(2)  A_4  D_4(\sqrt{-3}) $
& $2^{14}\cdot 3^7\cdot 5^2$      & 288   \\
& $y_1+\sqrt{-3}y_2+y_3+\frac{1}{60}x \ra$,
&
&
&                                  \\
& $y_1=\frac{1}{2}(1,1,1,1,1)$, $y_2=\frac{1}{5}(-4,1,1,1,1)$,
&
&
&                                         \\
& $y_3=\frac{1}{\sqrt{-3}}(1,1,1,1)$, $(x,x)=60$
&
&
&                                         \\ \hline
30
& $\la L^{(2)} \oplus \la x \ra \oplus \la z \ra,$
& $D_4(\sqrt{-3})^2  D_4(2)$
& $2^{14}\cdot 3^{11}$      & 288   \\
& $y_1 + \sqrt{-3}y_2+\frac{1}{3}x$,
$\frac{1}{2}(\sqrt{-3},1,1,1)+\frac{1}{2}x$,
&
&
&                                  \\
&
$\sqrt{-3}y_1 + y_2+\frac{1}{3}z$,
$\frac{1}{2}(-\sqrt{-3},1,1,1)+\frac{1}{2}z \ra$,
&
&
&                                         \\
& $y_{1,2}=\frac{1}{\sqrt{-3}}(1,1,1,1)$, $(x,x)=6$, $(z,z)=6$
&
&
&                                         \\ \hline
31
& $\la L^{(2)} \oplus \la x\ra \oplus \la z\ra ,$
& ${A_6}^2 $
& $2^{10}\cdot 3^5\cdot 5^2\cdot 7^2$      & 252   \\
& $y_1+\sqrt{-3}y_2+\frac{2}{7}x, \sqrt{-3}y_1+y_2+\frac{2}{7}z \ra$,
&
&
&                                  \\
& $y_i=\frac{1}{7}(-6,1,\dots,1)$, $(x,x)=7$, $(z,z)=7$
&
&
&                                         \\ \hline
32
& $\la L^{(2)} \oplus \la x \ra \oplus \la z \ra,$
& $A_6 A_2 D_4(2)$
& $2^{12}\cdot 3^6\cdot 5\cdot 7$      & 216   \\
& $3\sqrt{-3}y_1+y_2+\frac{1}{2}(\sqrt{-3},1,1,1)+\frac{1}{42}x+\frac{1}{3}z$,
&
&
&                                  \\
& $\frac{1}{2}(-\sqrt{-3},1,1,1)+\frac{\sqrt{-3}}{3}x-\frac{\sqrt{-3}}{6}z\ra$,
&
&
&                                  \\
& $y_1=\frac{1}{7}(-6,1,\dots,1)$, $y_3=\frac{1}{3}(-2,1,1)$
&
&
&                                         \\
& $(x,x)=42$, $(z,z)=6$
&
&
&                                         \\ \hline
33
& $\la L^{(2)},$
& ${A_5}^2 {A_2}^2 $
& $2^{12}\cdot 3^7\cdot 5^2$      & 216   \\
& $y_1+y_2+y_3+y_4, 2y_2+\omega y_3-\sqrt{-3}\omega y_4\ra$
&
&
&                                  \\
& $y_{1,2}=\frac{1}{6}(-5,1,\dots,1), y_{3,4}=\frac{1}{3}(-2,1,1)$,
&
&
&                                         \\ \hline
34
& $\la L^{(2)} \oplus \la x \ra ,$
& ${A_5}^2 A_3$
& $2^{13}\cdot 3^6\cdot 5^2$      & 216   \\
& $y_1+y_2+\frac{1}{3}x, y_3+\frac{\sqrt{-3}}{4}x, 2\sqrt{-3}(y_1-y_2)\ra$
&
&
&                                  \\
& $y_{1,2}=\frac{1}{6}(-5,1,\dots,1)$, $y_3=\frac{1}{4}(-3,1,1,1)$,
&
&
&                                         \\
& $(x,x)=12$
&
&
&                                         \\ \hline
35
& Complex conjugate of 34
&     &                    &          \\ \hline
36
& $\la L^{(2)} \oplus \la x\ra \oplus   \la z\ra ,$
& ${A_4}^2  D_4(2)$
& $2^{14}\cdot 3^4\cdot 5^2$      & 192   \\
& $2\sqrt{-3}y_1+y_2+\frac{2}{5}x, \frac{1}{2}(\sqrt{-3},1,1,1)+\frac{1}{2}x,$
&
&
&                                  \\
& $y_1+2\sqrt{-3}y_2+\frac{2}{5}z, \frac{1}{2}(1,-1,1,\sqrt{-3})+
\frac{1}{2}z \ra$,
&
&
&                                  \\
& $y_{1,2}=\frac{1}{5}(-4,1,\dots,1)$,
&
&
&                                  \\
& $(x,x)=10$, $(z,z)=10$
&
&
&                                         \\ \hline
\end{tabular}
\end{center}
\end{table}

\renewcommand{\arraystretch}{1.3}
\begin{table}[t]{Table 2. Unimodular lattices of rank 14 (continued)}
\label{table4}
\begin{center}

\begin{tabular}{|l|l|l|l|l|}  \hline
No
& $L$
& $L^{(2)}$       &  $|G(L)|$  &  $\mu_2(L)$     \\ \hline
37
& $\la L^{(2)}\oplus \la x\ra ,$
& $A_5  A_4  A_3  A_1$
& $2^{12}\cdot 3^5\cdot 5^2$      & 192   \\
& $2y_1 + \sqrt{-3}y_2 + y_3 + y_4 + \frac{1}{60}x$,
&
&
&                                         \\
& $3\omega y_1 + 2y_3 + y_4 \ra$, $(x,x)=60$,
&
&
&                                         \\
& $y_1=\frac{1}{6}(-5,1,\dots,1)$, $y_3=\frac{1}{4}(-3,1,1,1)$,
&
&
&                                         \\
& $y_2=\frac{1}{5}(-4,1,1,1,1)$, $y_4=\frac{1}{2}(-1,1)$
&
&
&                                         \\ \hline
38
& Complex conjugate of 37
&     &                    &          \\ \hline
39
& $$
& ${A_5}^2 {A_1}^2$
& $2^{12}\cdot 3^5\cdot 5^2$      & 192   \\ \hline
40
& $$
& ${D_3(\sqrt{-3})}^2 {A_2}^4$
& $2^{10}\cdot 3^{12}$      & 180   \\ \hline
41
& $$
& $A_5 D_3(\sqrt{-3}) {A_2}^2 $
& $2^9\cdot 3^9\cdot 5$      & 180   \\  \hline
42
& $$
& ${A_4}^2  D_3(\sqrt{-3})  A_1$
& $2^{10}\cdot 3^6\cdot 5^2$      & 180   \\ \hline
43
& $$
& $D_4(2) {D_3(\sqrt{-3})}^2 $
& $2^{13}\cdot 3^{10}$      & 180   \\ \hline
44
& $\la L^{(2)} \oplus\la x\ra \oplus \la z\ra ,$
& ${A_4}^2  {A_2}^2$
& $2^{10}\cdot 3^5\cdot 5^2$      & 156   \\
& $y_1+2\sqrt{-3}y_2+\frac{1}{5}x, y_3+y_4+\frac{1}{3}x,$
&
&
&                                  \\
& $2\sqrt{-3}y_1+y_2+\frac{1}{15}z, y_3-y_4+\frac{1}{3}z \ra$,
&
&
&                                  \\
& $y_{1,2}=\frac{1}{5}(-4,1,\dots,1)$, $(x,x)=15$,
&
&
&                                  \\
& $y_{3,4}=\frac{1}{3}(-2,1,1)$, $(z,z)=15$
&
&
&                                         \\ \hline
45
& $\la L^{(2)} \oplus \la x\ra \oplus \la z\ra \oplus \la u\ra ,$
& ${A_4}^2  A_3$
& $2^{11}\cdot 3^4\cdot 5^2$      & 156   \\
& $2\sqrt{-3}(y_1-y_2) + y_3 +\frac{1}{2}u +\frac{1}{20}x$,
&
&
&                                         \\
& $2\sqrt{-3}(y_1+y_2)+2y_3+\frac{\sqrt{-3}}{4}u+\frac{1}{20}z \ra$,
&
&
&                                  \\
& $y_{1,2}=\frac{1}{5}(-4,1,\dots,1)$, $(x,x)=20$,
&
&
&                                  \\
& $y_3=\frac{1}{4}(-3,1,1,1)$, $(z,z)=20$, $(u,u)=4$
&
&
&                                  \\ \hline
46
& $$
& $A_4 {A_3}^2  {A_1}^2$
& $2^{13}\cdot 3^4\cdot 5$      & 144   \\ \hline
47
& $$
& ${A_3}^4$
& $2^{15}\cdot 3^5$      & 144   \\ \hline
48
& $$
& ${A_3}^2  {A_2}^2  {A_1}^2$
& $2^{12}\cdot 3^5$      & 120   \\ \hline
49
& $$
& ${A_3}^3  {A_1}^2$
& $2^{13}\cdot 3^5$      & 120   \\ \hline
50
& Complex conjugate of 49
&     &                    &          \\ \hline
51
& $$
& ${A_2}^6 $
& $2^9\cdot 3^9$      & 108   \\ \hline
52
& $$
& $A_3  {A_2}^4$
& $2^{10}\cdot 3^7$      & 108   \\ \hline
53
& $$
& ${A_3}^2  {A_2}^2$
& $2^{11}\cdot 3^6$      & 108   \\ \hline
54
& $$
& ${A_2}^4  {A_1}^2$
& $2^{10}\cdot 3^5$      & 84   \\ \hline
55
& $$
& ${A_2}^2  {A_1}^6$
& $2^{11}\cdot 3^5$      & 72   \\ \hline
56
& $$
& ${A_2}^4$
& $2^8\cdot 3^9$      & 72   \\ \hline
57
& $$
& ${A_1}^8$
& $2^{13}\cdot 3^2\cdot 7$      & 48   \\ \hline
58
& $$
& ${A_1}^6$
& $2^{13}\cdot 3^4\cdot 5$      & 36   \\ \hline
\end{tabular}
\end{center}
\end{table}

\renewcommand{\arraystretch}{1.0}
\begin{table}[t]{Table 3. Unimodular lattices of rank 15}
\label{table11}
\begin{center}

\begin{tabular}{|l|l|l|l|l|}  \hline
No
& $L^{(2)}$       &  $|G(L)|$  &  $\mu_2(L)$     \\ \hline
1
& $\varnothing$             & $2^9\cdot 3^7\cdot 5\cdot 7$ & 0  \\ \hline
2
& $\varnothing$             & $2^5\cdot 3^7$ & 0  \\ \hline
3
& $D_{15}(\sqrt{-3})$
& $2^{12}\cdot 3^{20}\cdot 5^3\cdot 7^2\cdot 11\cdot 13$ & 1890  \\ \hline
4
& $D_9(\sqrt{-3})D_6(\sqrt{-3})$
& $2^{12}\cdot 3^{19}\cdot 5^2\cdot 7$      & 918  \\ \hline
5
& $E_7 D_7(\sqrt{-3})$
& $2^{15}\cdot 3^{13}\cdot 5^2\cdot 7^2$ & 756     \\ \hline
6
& $A_{15}$
& $2^{16}\cdot 3^7\cdot 5^3\cdot 7^2\cdot 11\cdot 13$      & 720  \\ \hline
7
& $A_{14}$
& $2^{12}\cdot 3^7\cdot 5^3\cdot 7^2\cdot 11\cdot 13$      & 630  \\ \hline
8
& complex conjugate of 7
&       &   \\ \hline
9
& $D_9(2)D_5(\sqrt{-3})$
& $2^{19}\cdot 3^{10}\cdot 5^2\cdot 7$      & 612  \\ \hline
10
& $A_8D_7(\sqrt{-3})$
& $2^{12}\cdot 3^{13}\cdot 5^2\cdot 7^2$      & 594  \\ \hline
11
& $D_6(\sqrt{-3})D_6(\sqrt{-3})D_3(\sqrt{-3})$
& $2^{11}\cdot 3^{17}\cdot 5^2$      & 594  \\ \hline
12
& $A_9U_5$
& $2^{16}\cdot 3^{9}\cdot 5^3\cdot 7$      & 540  \\ \hline
13
& $A_5D_5(\sqrt{-3})U_5$
& $2^{15}\cdot 3^{12}\cdot 5^3$      & 540  \\ \hline
14
& $E_6U_5D_3(\sqrt{-3})$
& $2^{16}\cdot 3^{13}\cdot 5^2$      & 540  \\ \hline
15
& $A_{11}D_4(\sqrt{-3})$
& $2^{14}\cdot 3^{10}\cdot 5^2\cdot 7\cdot 11$      & 504  \\ \hline
16
& $A_8D_6(\sqrt{-3})$
& $2^{12}\cdot 3^{12}\cdot 5^2\cdot 7$      & 486  \\ \hline
17 & $D_7(2)E_6$ & $2^{18}\cdot 3^8\cdot 5^2\cdot 7$ & 468 \\
\hline
18
& $A_{11}D_3(\sqrt{-3})$
& $2^{12}\cdot 3^9\cdot 5^2\cdot 7\cdot 11$      & 450  \\ \hline
19
& $A_9D_5(\sqrt{-3})$
& $2^{12}\cdot 3^{10}\cdot 5^3\cdot 7$      & 450  \\ \hline
20
& complex conjugate of 20
&       &   \\ \hline
21
& $A_8E_6$
& $2^{15}\cdot 3^9\cdot 5^2\cdot 7$      & 432  \\ \hline
22
& $D_6(\sqrt{-3})D_3^3(\sqrt{-3})$
& $2^9\cdot 3^{18}\cdot 5$      & 432  \\ \hline
23
& $U_5D_4(\sqrt{-3})D_4(\sqrt{-3})$
& $2^{15}\cdot 3^{13}\cdot 5$      & 432  \\ \hline
24
& $D_7(2)A_7$
& $2^{18}\cdot 3^5\cdot 5^2\cdot 7^2$      & 420  \\ \hline
25
& $D_6(2)D_5(\sqrt{-3})A_3$
& $2^{16}\cdot 3^9\cdot 5^2$      & 396  \\ \hline
26
& $A_{10}A_4$
& $2^{12}\cdot 3^6\cdot 5^3\cdot 7\cdot 11$      & 390  \\ \hline
27
& $A_5D_5(\sqrt{-3})D_4(\sqrt{-3})$
& $2^{11}\cdot 3^{12}\cdot 5^2$      & 378  \\ \hline
28
& $A_9A_5$
& $2^{13}\cdot 3^7\cdot 5^3\cdot 7$      & 360  \\ \hline
29
& complex conjugate of 28
&       &  \\ \hline
30
& $A_7D_6(2)$
& $2^{17}\cdot 3^5\cdot 5^2\cdot 7$      & 348  \\ \hline
31
& $A_6D_5(\sqrt{-3})A_3$
& $2^{11}\cdot 3^9\cdot 5^2\cdot 7$      & 342  \\ \hline
32
& $A_8 A_6$
& $2^{12}\cdot 3^7\cdot 5^2\cdot 7^2$      & 342  \\ \hline
33
& complex conjugate of 32
&       &  \\ \hline
34
& $A_8D_4(\sqrt{-3})A_2$
& $2^{12}\cdot 3^{10}\cdot 5\cdot 7$      & 342  \\ \hline
35
& complex conjugate of 34
&       &  \\ \hline
36
& $A_8D_5(2)$
& $2^{15}\cdot 3^6\cdot 5^2\cdot 7$      & 336  \\ \hline
37
& $D_5(\sqrt{-3})D_4(2)D_4(2)$
& $2^{17}\cdot 3^{10}\cdot 5$      & 324  \\ \hline
38
& $D_6(2)A_5D_3(\sqrt{-3})$
& $2^{15}\cdot 3^8\cdot 5^2$      & 324  \\ \hline
39
& $A_8D_3(\sqrt{-3})A_3$
& $2^{12}\cdot 3^9\cdot 5\cdot 7$      & 306  \\ \hline
40
& $A_8A_5$
& $2^{12}\cdot 3^7\cdot 5^2\cdot 7$      & 306  \\ \hline
41
& complex conjugate of 40
&       &  \\ \hline
42
& $A_7A_6$
& $2^{12}\cdot 3^5\cdot 5^2\cdot 7^2$      & 294  \\ \hline
43
& complex conjugate of 42
&       &  \\ \hline
\end{tabular}
\end{center}
\end{table}

\begin{table}[t]{Table 3. Unimodular lattices of rank 15 (continued)}
\label{table12}
\begin{center}

\begin{tabular}{|l|l|l|l|l|}  \hline
No
& $L^{(2)}$       &  $|G(L)|$  &  $\mu_2(L)$     \\ \hline
44
& $D_5(2)D_4(\sqrt{-3})A_4$
& $2^{14}\cdot 3^7\cdot 5^2$      & 288  \\ \hline
45
& $D_4(\sqrt{-3})A_5A_5$
& $2^{12}\cdot 3^9\cdot 5^2$      & 288  \\ \hline
46
& $D_5(2)D_5(2)A_3$
& $2^{19}\cdot 3^4\cdot 5^2$      & 276  \\ \hline
47
& $A_6A_5D_3(\sqrt{-3})$
& $2^{10}\cdot 3^8\cdot 5^2\cdot 7$      & 270  \\ \hline
48
& $A_6D_4(\sqrt{-3})A_3$
& $2^{11}\cdot 3^{8}\cdot 5\cdot 7$      & 270  \\ \hline
49
& $A_5D_4(\sqrt{-3})D_3(\sqrt{-3})A_2$
& $2^{10}\cdot 3^{11}\cdot 5$      & 270  \\ \hline
50
& $D_4(\sqrt{-3})D_4(\sqrt{-3})A_2A_2A_2$
& $2^{11}\cdot 3^{13}$             & 270  \\ \hline
51
& $D_3(\sqrt{-3})^5$
& $2^8\cdot 3^{16}\cdot 5$      & 270  \\ \hline
52
& $A_6A_6A_2$
& $2^{11}\cdot 3^6\cdot 5^2\cdot 7^2$      & 270  \\ \hline
53
& complex conjugate of 52
&       &  \\ \hline
54
& $A_7A_4A_3$
& $2^{14}\cdot 3^5\cdot 5^2\cdot 7$      & 264  \\ \hline
55
& $A_6D_5(2)A_2$
& $2^{13}\cdot 3^5\cdot 5^2\cdot 7$      & 264  \\ \hline
56
& complex conjugate of 55
&       &  \\ \hline
57
& $A_6A_6A_1$
& $2^{10}\cdot 3^5\cdot 5^2\cdot 7^2$      & 258  \\ \hline
58
& $D_5(2)D_4(2)D_3(\sqrt{-3})A_1$
& $2^{16}\cdot 3^7\cdot 5$           & 252  \\ \hline
59
& $A_5A_5D_4(2)$
& $2^{16}\cdot 3^7\cdot 5^2$      & 252  \\ \hline
60
& $A_7A_4A_2$
& $2^{12}\cdot 3^5\cdot 5^2\cdot 7$      & 246  \\ \hline
61
& $D_5(2)A_4A_4$
& $2^{14}\cdot 3^4\cdot 5^3$      & 240  \\ \hline
62
& $D_5(2)A_4A_4$
& $2^{15}\cdot 3^4\cdot 5^3$      & 240  \\ \hline
63
& $A_7A_3A_3$
& $2^{15}\cdot 3^5\cdot 5\cdot 7$      & 240  \\ \hline
64
& $A_5D_4(\sqrt{-3})A_2A_2$
& $2^{10}\cdot 3^9\cdot 5$        & 234  \\ \hline
65
& $A_4A_4D_4(\sqrt{-3})A_1$
& $2^{11}\cdot 3^7\cdot 5^2$      & 234  \\ \hline
66
& $A_6A_5A_2$
& $2^{10}\cdot 3^6\cdot 5^2\cdot 7$      & 234  \\ \hline
67
& $A_6A_5A_2$
& $2^{10}\cdot 3^6\cdot 5^2\cdot 7$      & 234  \\ \hline
68
& complex conjugate of 67
&       &  \\ \hline
69
& $A_5A_5D_3(\sqrt{-3})$
& $2^{11}\cdot 3^8\cdot 5^2$      & 234  \\ \hline
70
& complex conjugate of 69
&       &  \\ \hline
71
& $D_5(2)A_3A_3A_3$
& $2^{17}\cdot 3^6\cdot 5$      & 228  \\ \hline
72
& $A_6A_4A_3$
& $2^{11}\cdot 3^5\cdot 5^2\cdot 7$      & 222  \\ \hline
73
& complex conjugate of 72
&       &  \\ \hline
74
& $D_4(\sqrt{-3})A_3A_3A_3$
& $2^{14}\cdot 3^9$            & 216  \\ \hline
75
& $D_4(2)D_4(\sqrt{-3})A_2A_2$
& $2^{13}\cdot 3^{10}$            & 216  \\ \hline
76
& $D_4(\sqrt{-3})D_3(\sqrt{-3})A_2A_2A_2$
& $2^9\cdot 3^{13}$            & 216  \\ \hline
77
& $A_5D_3(\sqrt{-3})^2A_2$
& $2^9\cdot 3^{11}\cdot 5$            & 216  \\ \hline
78
& $A_6D_3(\sqrt{-3})A_2A_2$
& $2^9\cdot 3^8\cdot 5\cdot 7$            & 216  \\ \hline
79
& $A_6D_4(2)A_2$
& $2^{12}\cdot 3^6\cdot 5\cdot 7$            & 216  \\ \hline
80
& $A_5A_5A_3$
& $2^{12}\cdot 3^6\cdot 5^2$            & 216  \\ \hline
81
& $A_5A_5A_3$
& $2^{13}\cdot 3^6\cdot 5^2$            & 216  \\ \hline
82
& complex conjugate of 81
&       &  \\ \hline
83
& $A_6A_4A_2A_1$
& $2^{10}\cdot 3^5\cdot 5^2\cdot 7$            & 210  \\ \hline
84
& complex conjugate of 83
&       &  \\ \hline
85
& $A_5A_4A_4$
& $2^{11}\cdot 3^5\cdot 5^3$            & 210  \\ \hline
86
& complex conjugate of 85
&       &  \\ \hline
87
& $A_5D_4(2)A_3A_1$
& $2^{15}\cdot 3^5\cdot 5$            & 204  \\ \hline
\end{tabular}
\end{center}
\end{table}

\begin{table}[t]{Table 3. Unimodular lattices of rank 15 (continued)}
\label{table13}
\begin{center}
\begin{tabular}{|l|l|l|l|l|}  \hline
No
& $L^{(2)}$       &  $|G(L)|$  &  $\mu_2(L)$     \\ \hline
88
& $A_6A_3A_2A_2$
& $2^{10}\cdot 3^6\cdot 5\cdot 7$      & 198  \\ \hline
89
& $A_5A_5A_2$
& $2^{11}\cdot 3^7\cdot 5^2$      & 198  \\ \hline
90
& $A_5A_5A_2$
& $2^{11}\cdot 3^6\cdot 5^2$      & 198  \\ \hline
91
& complex conjugate of 90
&       &  \\ \hline
92
& $A_5A_3D_3(\sqrt{-3})A_2$
& $2^{10}\cdot 3^8\cdot 5$      & 198  \\ \hline
93
& complex conjugate of 92
&       &  \\ \hline
94
& $A_5A_4A_3A_1$
& $2^{12}\cdot 3^5\cdot 5^2$      & 192  \\ \hline
95
& complex conjugate of 94
&       &  \\ \hline
96
& $A_4A_4D_4(2)$
& $2^{14}\cdot 3^4\cdot 5^2$      & 192  \\ \hline
97
& complex conjugate of 96
&       &  \\ \hline
98
& $A_5A_4A_3$
& $2^{11}\cdot 3^5\cdot 5^2$      & 186  \\ \hline
99
& $A_5A_4A_3$
& $2^{11}\cdot 3^5\cdot 5^2$      & 186  \\ \hline
100
& $A_5D_4(2)A_1A_1A_1$
& $2^{14}\cdot 3^6\cdot 5$      & 180  \\ \hline
101
& $D_4(2)D_3(\sqrt{-3})A_3A_1A_1A_1$
& $2^{14}\cdot 3^7$             & 180  \\ \hline
102
& $D_4D_4(2)A_3$
& $2^{17}\cdot 3^6$         & 180  \\ \hline
103
& $A_3A_3A_3A_3A_3$
& $2^{18}\cdot 3^7\cdot 5$      & 180  \\ \hline
104
& $A_4A_4A_4$
& $2^{10}\cdot 3^5\cdot 5^3$      & 180  \\ \hline
105
& complex conjugate of 104
&       &  \\ \hline
106
& $A_5A_4A_2A_1$
& $2^{10}\cdot 3^5\cdot 5^2$      & 174  \\ \hline
107
& $A_5A_4A_2A_1$
& $2^{10}\cdot 3^5\cdot 5^2$      & 174  \\ \hline
108
& complex conjugate of 107
&       &  \\ \hline
109
& $A_4A_4A_3A_2$
& $2^{11}\cdot 3^5\cdot 5^2$      & 174  \\ \hline
110
& complex conjugate of 109
&       &  \\ \hline
111
& $A_4D_4(2)A_2A_2$
& $2^{12}\cdot 3^5\cdot 5$           & 168  \\ \hline
112
& $A_4D_4(2)A_2A_2$
& $2^{13}\cdot 3^5\cdot 5$           & 168  \\ \hline
113
& $A_4A_3A_3A_3$
& $2^{14}\cdot 3^5\cdot 5$        & 168  \\ \hline
114
& $A_5A_3A_3A_1$
& $2^{13}\cdot 3^5\cdot 5$      & 168  \\ \hline
115
& complex conjugate of 114
&       &  \\ \hline
116
& $A_4D_3(\sqrt{-3})A_3A_1A_1$
& $2^{10}\cdot 3^6\cdot 5$        & 162  \\ \hline
117
& $A_3A_3A_3D_3(\sqrt{-3})$
& $2^{11}\cdot 3^8$            & 162  \\ \hline
118
& $A_3A_3D_3(\sqrt{-3})A_2A_2$
& $2^{11}\cdot 3^8$            & 162  \\ \hline
119
& $A_5A_2A_2A_2A_2$
& $2^9\cdot 3^8\cdot 5$        & 162  \\ \hline
120
& $D_3(\sqrt{-3})D_3(\sqrt{-3})A_2A_2A_2$
& $2^7\cdot 3^{12}$            & 162  \\ \hline
121
& $A_5A_3A_2A_2$
& $2^{10}\cdot 3^6\cdot 5$            & 162  \\ \hline
122
& complex conjugate of 121
&       &  \\ \hline
123
& $A_4A_4A_3A_1$
& $2^{11}\cdot 3^4\cdot 5^2$      & 162  \\ \hline
124
& complex conjugate of 123
&       &  \\ \hline
125
& $A_4A_4A_2A_2$
& $2^{10}\cdot 3^5\cdot 5^2$            & 156  \\ \hline
126
& $D_4(2)A_3A_3A_1A_1$
& $2^{16}\cdot 3^4$               & 156  \\ \hline
127
& $A_4A_4A_2A_1A_1$
& $2^{10}\cdot 3^4\cdot 5^2$            & 150  \\ \hline
128
& $A_4A_4A_2A_1A_1$
& $2^{11}\cdot 3^4\cdot 5^2$            & 150  \\ \hline
129
& $A_4A_3A_3A_2$
& $2^{11}\cdot 3^5\cdot 5$            & 150  \\ \hline
130
& $A_4A_3A_3A_2$
& $2^{11}\cdot 3^5\cdot 5$            & 150  \\ \hline
131
& $A_4A_3A_3A_2$
& $2^{11}\cdot 3^5\cdot 5$            & 150  \\ \hline
132
& complex conjugate of 131
&       &  \\ \hline
\end{tabular}
\end{center}
\end{table}

\begin{table}[t]{Table 3. Unimodular lattices of rank 15 (continued)}
\label{table14}
\begin{center}

\begin{tabular}{|l|l|l|l|l|}  \hline
No
& $L^{(2)}$       &  $|G(L)|$  &  $\mu_2(L)$     \\ \hline
133
& $A_4A_3A_3A_1A_1$
& $2^{12}\cdot 3^4\cdot 5$      & 144  \\ \hline
134
& $A_3D_3(\sqrt{-3})A_3A_2A_2A_2$
& $2^8\cdot 3^9$      & 144  \\ \hline
135
& $A_5A_2A_2A_1A_1A_1$
& $2^{11}\cdot 3^6\cdot 5$      & 144  \\ \hline
136
& $A_4D_3(\sqrt{-3})A_2A_1A_1$
& $2^9\cdot 3^7\cdot 5$      & 144  \\ \hline
137
& $A_5A_2A_2A_2$
& $2^9\cdot 3^8\cdot 5$      & 144  \\ \hline
138
& $A_3A_3A_3A_3$
& $2^{16}\cdot 3^6$      & 144  \\ \hline
139
& $D_4(2)A_2A_2A_2A_2$
& $2^{12}\cdot 3^8$      & 144  \\ \hline
140
& $D_4(2)A_2A_2A_2A_2$
& $2^{13}\cdot 3^7$      & 144  \\ \hline
141
& $A_3A_3A_3A_2A_2$
& $2^{13}\cdot 3^7$      & 144  \\ \hline
142
& complex conjugate of 141
&       &  \\ \hline
143
& $A_4A_4A_1A_1A_1$
& $2^{11}\cdot 3^3\cdot 5^2$      & 138  \\ \hline
144
& $A_4A_3A_2A_2A_1$
& $2^{10}\cdot 3^5\cdot 5$      & 138  \\ \hline
145
& $A_4A_3A_2A_2A_1$
& $2^{10}\cdot 3^5\cdot 5$      & 138  \\ \hline
146
& complex conjugate of 145
&       &  \\ \hline
147
& $A_4A_3A_2A_2A_1$
& $2^{10}\cdot 3^5\cdot 5$      & 138  \\ \hline
148
& complex conjugate of 147
&       &  \\ \hline
149
& $A_4A_3A_3A_1$
& $2^{11}\cdot 3^4\cdot 5$      & 138  \\ \hline
150
& complex conjugate of 149
&       &  \\ \hline
151
& $A_4A_3A_2A_2$
& $2^9\cdot 3^5\cdot 5$      & 132  \\ \hline
152
& $A_3D_3(\sqrt{-3})A_2A_1A_1A_1$
& $2^9\cdot 3^7$             & 126  \\ \hline
153
& $D_3(\sqrt{-3})A_2A_2A_2A_2$
& $2^8\cdot 3^9$             & 126  \\ \hline
154
& $A_3A_3A_3A_2$
& $2^{11}\cdot 3^6$             & 126  \\ \hline
155
& $A_3A_3A_2A_2A_2$
& $2^{11}\cdot 3^6$         & 126  \\ \hline
156
& $A_3A_3A_2A_2A_2$
& $2^{11}\cdot 3^6$         & 126  \\ \hline
157
& complex conjugate of 156
&       &  \\ \hline
158
& $A_3A_3A_2A_2A_2$
& $2^{11}\cdot 3^7$         & 126  \\ \hline
159
& complex conjugate of 158
&       &  \\ \hline
160
& $A_4A_3A_2A_1A_1$
& $2^{10}\cdot 3^4\cdot 5$         & 126  \\ \hline
161
& complex conjugate of 160
&       &  \\ \hline
162
& $A_3A_3A_2A_2A_1A_1$
& $2^{11}\cdot 3^5$              & 120  \\ \hline
163
& $A_4A_3A_1A_1A_1A_1$
& $2^{12}\cdot 3^3\cdot 5$              & 120  \\ \hline
164
& $A_4A_3A_3A_3A_1$
& $2^8\cdot 3^6\cdot 5$              & 120  \\ \hline
165
& $A_3A_3A_3A_1A_1$
& $2^{12}\cdot 3^5$              & 120  \\ \hline
166
& $A_3A_3A_3A_1A_1$
& $2^{13}\cdot 3^4$              & 120  \\ \hline
167
& complex conjugate of 166
&       &  \\ \hline
168
& $A_4A_2A_2A_1A_1A_1$
& $2^9\cdot 3^4\cdot 5$        & 114  \\ \hline
169
& $A_3A_3A_2A_2A_1$
& $2^{10}\cdot 3^5$         & 114  \\ \hline
170
& $A_3A_3A_2A_2A_1$
& $2^{10}\cdot 3^5$         & 114  \\ \hline
171
& complex conjugate of 170
&       &  \\ \hline
172
& $A_3A_3A_2A_2A_1$
& $2^{11}\cdot 3^5$        & 114  \\ \hline
173
& complex conjugate of 172
&       &  \\ \hline
\end{tabular}
\end{center}
\end{table}

\begin{table}[t]{Table 3. Unimodular lattices of rank 15 (continued)}
\label{table15}
\begin{center}

\begin{tabular}{|l|l|l|l|l|}  \hline
No
& $L^{(2)}$       &  $|G(L)|$  &  $\mu_2(L)$     \\ \hline
174
& $A_4A_2A_2A_1A_1$
& $2^9\cdot 3^5\cdot 5$      & 108  \\ \hline
175
& $A_3A_3A_2A_2$
& $2^{10}\cdot 3^6$      & 108  \\ \hline
176
& $A_3A_2A_2A_2A_2$
& $2^9\cdot 3^6$      & 108  \\ \hline
177
& $A_3A_3A_1A_1A_1A_1A_1A_1$
& $2^{14}\cdot 3^4$      & 108  \\ \hline
178
& $D_3(\sqrt{-3})A_2A_2A_2$
& $2^6\cdot 3^{11}$      & 108  \\ \hline
179
& $D_3(\sqrt{-3})D_3(\sqrt{-3})$
& $2^{10}\cdot 3^{13}$      & 108  \\ \hline
180
& $D_3(\sqrt{-3})A_1A_1A_1A_1A_1A_1A_1A_1A_1$
& $2^{14}\cdot 3^7$      & 108  \\ \hline
181
& $A_2A_2A_2A_2A_2A_2$
& $2^8\cdot 3^9$      & 108  \\ \hline
182
& $A_2A_2A_2A_2A_2A_2$
& $2^{10}\cdot 3^{10}$      & 108  \\ \hline
183
& $A_3A_3A_2A_1A_1$
& $2^{10}\cdot 3^4$      & 102  \\ \hline
184
& $A_3A_2A_2A_2A_1A_1$
& $2^9\cdot 3^5$      & 102  \\ \hline
185
& $A_3A_2A_2A_2A_1A_1$
& $2^9\cdot 3^5$      & 102  \\ \hline
186
& $A_3A_2A_2A_2A_1A_1$
& $2^9\cdot 3^5$      & 102  \\ \hline
187
& complex conjugate of 186
&       &  \\ \hline
188
& $A_3A_3A_2A_1A_1$
& $2^{11}\cdot 3^4$      & 102  \\ \hline
189
& complex conjugate of 188
&       &  \\ \hline
190
& $A_3A_3A_2A_1A_1$
& $2^{11}\cdot 3^4$      & 102  \\ \hline
191
& complex conjugate of 190
&       &  \\ \hline
192
& $A_3A_2A_2A_1A_1A_1A_1$
& $2^{11}\cdot 3^4$      & 96  \\ \hline
193
& $A_3A_3A_1A_1A_1A_1$
& $2^{13}\cdot 3^3$      & 96  \\ \hline
194
& $A_3A_2A_2A_2A_1$
& $2^9\cdot 3^5$      & 96  \\ \hline
195
& $A_3A_2A_2A_1A_1A_1A_1$
& $2^{12}\cdot 3^4$      & 96  \\ \hline
196
& $A_3A_3A_1A_1A_1$
& $2^{11}\cdot 3^4$      & 90  \\ \hline
197
& $A_2A_2A_2A_2A_1A_1A_1$
& $2^9\cdot 3^6$      & 90  \\ \hline
198
& $A_3A_2A_2A_1A_1A_1$
& $2^9\cdot 3^5$      & 90  \\ \hline
199
& $A_3A_2A_2A_1A_1A_1$
& $2^9\cdot 3^4$      & 90  \\ \hline
200
& $A_3A_2A_2A_1A_1A_1$
& $2^9\cdot 3^4$      & 90  \\ \hline
201
& complex conjugate of 200
&       &  \\ \hline
202
& $A_2A_2A_2A_2A_2$
& $2^8\cdot 3^7$      & 90  \\ \hline
203
& $A_2A_2A_2A_2A_2$
& $2^7\cdot 3^7$      & 90  \\ \hline
204
& complex conjugate of 203
&       &  \\ \hline
205
& $A_2A_2A_2A_2A_2$
& $2^8\cdot 3^7\cdot 5$      & 90  \\ \hline
206
& complex conjugate of 205
&       &  \\ \hline
207
& $A_3A_2A_2A_1A_1$
& $2^9\cdot 3^4$      & 84  \\ \hline
208
& $A_3A_1A_1A_1A_1A_1A_1A_1A_1$
& $2^{15}\cdot 3^3$      & 84  \\ \hline
209
& $A_2A_2A_2A_2A_1A_1$
& $2^8\cdot 3^5$      & 84  \\ \hline
210
& complex conjugate of 209
&       &  \\ \hline
211
& $A_3A_2A_1A_1A_1A_1$
& $2^9\cdot 3^3$      & 78  \\ \hline
212
& $A_2A_2A_2A_2A_1$
& $2^7\cdot 3^5$      & 78  \\ \hline
213
& $A_2A_2A_2A_1A_1A_1A_1$
& $2^9\cdot 3^4$      & 78  \\ \hline
214
& $A_2A_2A_2A_1A_1A_1A_1$
& $2^9\cdot 3^4$      & 78  \\ \hline
215
& complex conjugate of 214
&       &  \\ \hline
\end{tabular}
\end{center}
\end{table}

\begin{table}[t]{Table 3. Unimodular lattices of rank 15 (continued)}
\label{table16}
\begin{center}

\begin{tabular}{|l|l|l|l|l|}  \hline
No
& $L^{(2)}$       &  $|G(L)|$  &  $\mu_2(L)$     \\ \hline
216
& $A_3A_2A_2$
& $2^8\cdot 3^6$      & 72  \\ \hline
217
& $A_3A_2A_1A_1A_1$
& $2^9\cdot 3^5$      & 72  \\ \hline
218
& $A_2A_2A_2A_1A_1A_1$
& $2^7\cdot 3^5$      & 72  \\ \hline
219
& $A_3A_3A_3A_1A_1A_1$
& $2^8\cdot 3^6$      & 72  \\ \hline
220
& $A_2A_2A_1A_1A_1A_1A_1A_1$
& $2^{10}\cdot 3^5$      & 72  \\ \hline
221
& $A_2A_2A_1A_1A_1A_1A_1A_1$
& $2^{10}\cdot 3^4$      & 72  \\ \hline
222
& complex conjugate of 221
&       &  \\ \hline
223
& $A_3A_1A_1A_1A_1A_1A_1$
& $2^{12}\cdot 3^3$      & 72  \\ \hline
224
& complex conjugate of 223
&       &  \\ \hline
225
& $A_2A_2A_2A_1A_1$
& $2^7\cdot 3^4$      & 66  \\ \hline
226
& $A_2A_2A_2A_1A_1$
& $2^7\cdot 3^5$      & 66  \\ \hline
227
& $A_2A_2A_1A_1A_1A_1A_1$
& $2^8\cdot 3^3$      & 66  \\ \hline
228
& $A_2A_2A_1A_1A_1A_1A_1$
& $2^9\cdot 3^3$      & 66  \\ \hline
229
& complex conjugate of 228
&       &  \\ \hline
230
& $A_2A_2A_1A_1A_1A_1$
& $2^8\cdot 3^3$      & 60  \\ \hline
231
& $A_2A_2A_1A_1A_1A_1$
& $2^{10}\cdot 3^3$      & 60  \\ \hline
232
& $A_1A_1A_1A_1A_1A_1A_1A_1A_1A_1$
& $2^{13}\cdot 3^2\cdot 5$      & 60  \\ \hline
233
& $A_2A_2A_2$
& $2^5\cdot 3^8$      & 54  \\ \hline
234
& $A_2A_2A_1A_1A_1$
& $2^7\cdot 3^5$      & 54  \\ \hline
235
& $A_2A_2A_1A_1A_1$
& $2^7\cdot 3^4$      & 54  \\ \hline
236
& complex conjugate of 235
&       &  \\ \hline
237
& $A_2A_1A_1A_1A_1A_1A_1$
& $2^9\cdot 3^2$      & 54  \\ \hline
238
& $A_2A_1A_1A_1A_1A_1A_1$
& $2^9\cdot 3^3$      & 54  \\ \hline
239
& complex conjugate of 238
&       &  \\ \hline
240
& $A_2A_1A_1A_1A_1A_1$
& $2^8\cdot 3^3$      & 48  \\ \hline
241
& $A_1A_1A_1A_1A_1A_1A_1A_1$
& $2^{14}\cdot 3$      & 48  \\ \hline
242
& $A_1A_1A_1A_1A_1A_1A_1A_1$
& $2^{10}\cdot 3^2\cdot 7$      & 48  \\ \hline
243
& complex conjugate of 242
&       &  \\ \hline
244
& $A_2A_1A_1A_1A_1$
& $2^8\cdot 3^2$      & 42  \\ \hline
245
& $A_1A_1A_1A_1A_1A_1A_1$
& $2^9\cdot 3^2$      & 42  \\ \hline
246
& $A_1A_1A_1A_1A_1A_1A_1$
& $2^9\cdot 3\cdot 7$      & 42  \\ \hline
247
& complex conjugate of 246
&       &  \\ \hline
248
& $A_2A_1A_1A_1$
& $2^6\cdot 3^5$      & 36  \\ \hline
249
& $A_1A_1A_1A_1A_1A_1$
& $2^9\cdot 3^2$      & 36  \\ \hline
250
& $A_2A_2$
& $2^6\cdot 3^7$      & 36  \\ \hline
251
& $A_3$
& $2^{13}\cdot 3^6$      & 36  \\ \hline
252
& $A_1A_1A_1A_1A_1$
& $2^7\cdot 3\cdot 5$      & 30  \\ \hline
253
& $A_1A_1A_1A_1A_1$
& $2^8\cdot 3^2\cdot 5$      & 30  \\ \hline
254
& complex conjugate of 253
&       &  \\ \hline
255
& $A_1A_1A_1A_1$
& $2^9\cdot 3^2$      & 24  \\ \hline
256
& $A_1A_1A_1$
& $2^5\cdot 3^4$      & 18  \\ \hline
257
& $A_2$
& $2^4\cdot 3^6\cdot 7\cdot 13$      & 18  \\ \hline
258
& $A_1A_1$
& $2^{10}\cdot 3^3$      & 12  \\ \hline
259
& $A_1$
& $2^4\cdot 3^2\cdot 7\cdot 13$      & 6  \\ \hline
\end{tabular}
\end{center}
\end{table}


\end{document}